\documentclass[11pt,a4paper,oneside,reqno]{amsart}

\usepackage[all]{xy}
\SelectTips{cm}{}  \calclayout
\usepackage{amsfonts,amssymb,amscd,amsmath,enumerate,verbatim,calc}

\usepackage{amscd,amssymb,amsopn,amsmath,amsthm,graphics,amsfonts,enumerate,verbatim,calc}
\usepackage[dvips]{graphicx}
\usepackage[colorlinks=true,linkcolor=blue,citecolor=blue]{hyperref}
\input xy
\xyoption{all}

\newcommand{\CA}{\mathcal{A} }

\newcommand{\CF}{\mathcal{F} }
\newcommand{\CG}{\mathcal{G} }

\newcommand{\CH}{\mathcal{H} }

\newcommand{\CK}{\mathcal{K} }

\newcommand{\Hom}{{\rm{Hom}}}

\newtheorem{theorem}{Theorem}[section]
\newtheorem{corollary}[theorem]{Corollary}
\newtheorem{lemma}[theorem]{Lemma}

\theoremstyle{definition}
\newtheorem{definition}[theorem]{Definition}

\theoremstyle{plain}

\theoremstyle{definition}

\numberwithin{equation}{section}
\begin{document}

\title[Multiplicative (generalized)-derivations of prime rings]{Multiplicative (generalized)-derivations of prime rings that act as $n-$(anti)homomorphisms}
\author[G. S. Sandhu]{G. S. Sandhu}

\address{Department of Mathematics, Patel Memorial National College, Rajpura 140401, Punjab, India.}
\email{gurninder\_rs@pbi.ac.in}

\begin{abstract}
Let $R$ be a prime ring. In this note, we describe the possible forms of multiplicative (generalized)-derivations of $R$
that act as $n-$homomorphism or $n-$antihomomorphism on nonzero ideals of $R.$ Consequently, from the given results one can
easily deduce the results of Gusi\'{c} \cite{Gusic2005}.
\end{abstract}
\maketitle

\section{Introduction}
Throughout this paper, $R$ will always denote an associative prime ring with center $Z(R)$ and $C$ the extended centroid of $R.$ It is well-known that in this case $C$ is a field. For any $x,y\in R,$ the symbol $[x,y]$ denotes the commutator $xy-yx.$ Recall, a ring is said to be prime if $xRy=(0)$ (where $x,y\in R$) implies $x=0$ or $y=0.$ An additive mapping $d:R\to R$ is said to be a derivation if $d(xy)=d(x)y+xd(y)$ for all $x,y\in R.$ In 1991, Bre$\check{s}$ar \cite{Bresar1991} introduced the notion of generalized derivation as follows: an additive mapping $F:R\to R$ is said to be a generalized derivation if $F(xy)=F(x)y+xd(y)$ for all $x,y\in R,$ where $d$ is a derivation of $R.$ The concept of generalized derivation covers both the notions of derivation and left multiplier (i.e., an additive mapping $T:R\to R$ satisfying $T(xy)=T(x)y$ for all $x,y\in R$). Now if we relax the assumption of additivity in the notion of derivation, then it is called multiplicative derivation, i.e., a mapping $\delta:R\to R$ (not necessarily additive) satisfying $\delta(xy)=\delta(x)y+x\delta(y)$ for all $x,y\in R.$ Recently, Dhara and Ali \cite{Dhara2013} extended the notion of multiplicative derivation to multiplicative (generalized)-derivation. Accordingly, a mapping $F:R\to R$ (not necessarily additive) is said to be a multiplicative (generalized)-derivation of $R$ if $F(xy)=F(x)y+x\delta(y)$ for all $x,y\in R,$ where $\delta$ is a multiplicative derivation of $R.$ Clearly, every generalized derivation is a multiplicative (generalized)-derivation, however the converse is not generally true ( see \cite{Dhara2013}, Example 1.1). Recall that a mapping $f$ of $R$ is said to be acts as a homomorphism (resp. anti-homomorphism) on an appropriate subset $K$ of $R$ if $f(xy)=f(x)f(y)$ (resp. $f(xy)=f(y)f(x)$) for all $x,y\in K.$ Following Hezajian et al. \cite{Hejazian2005}, a mapping $f$ of $R$ is said to be acts a an $n-$homomorphism (resp. $n-$antihomomorphism) of $R$ if for any $x_{i}\in R,$ where $i=1,2,\cdots,n$; $f(\prod_{i=1}^{n}x_{i})=\prod_{i=1}^{n}f(x_{i})$ (resp. $f(\prod_{i=1}^{n}x_{i})=f(x_{n})f(x_{n-1})\cdots f(x_{1})$). Initially, the notion of $n-$homomorphism was introduced and studied for complex algebras by Hejazian et al. \cite{Hejazian2005}, where some significant properties of $n-$homomorphisms are discussed on Banach algebras. Moreover, it is not difficult to see that every homomorphism of $R$ is $n-$homomorphism (for $n>  2$), but the converse is not necessarily true (see \cite{Hejazian2005}).

\par Till date, there exist many results in the literature showing that the global structure of $R$ is often tightly connected to the behaviour of additive mappings defined on $R.$ In 1989, a result due to Bell and Kappe \cite{Bell1989} states that if a prime ring $R$ admits a derivation $d$ that acts as homomorphism or anti-homomorphism on a nonzero right ideal $U$ of $R,$ then $d=0.$ Later Asma et al. \cite{Asma2003} proved that this result also holds on nonzero square-closed Lie ideals of prime rings. Moreover, Rehman \cite{Rehman2004} established this result for generalized derivations of prime rings. In fact, he proved that if $F$ is a nonzero generalized derivation of a 2-torsion free prime ring $R$ that acts as homomorphism or anti-homomorphism on a nonzero ideal of $R$ and $d\neq 0,$ then $R$ is commutative. Recently, Lukashenko \cite{Lukashenko2015} provided a new direction to these studies by investigating derivations acting as homomorphisms or anti-homomorphisms in differentially semiprime rings. Now it seems interesting to extend the results of generalized derivations to multiplicative (generalized)-derivations. In this context, Gusi\'{c} \cite{Gusic2005} gave the complete form of Rehman's result as follows: \emph{Let $R$ be an associative prime ring, $F$ be a multiplicative (generalized)-derivation of $R$ associated with a multiplicative derivation $\delta$ and $I$ be a nonzero ideal of $R.$
\begin{itemize}
\item[(a)] Assume that $F$ acts as homomorphism on $I.$ Then $\delta=0,$ and $F=0$ or $F(x)=x$ for all $x\in R.$
\item[(b)] Assume that $F$ acts as anti-homomorphism on $I.$ Then $\delta=0,$ and $F=0$ or $F(x)=x$ for all $x\in R$ (in this case $R$ should be commutative).
\end{itemize}
}
 In view of our above discussion, we find it reasonable to extend the results of derivations acting as homomorphisms (resp. anti-homomorphisms) to $n-$homomorphisms (resp. $n-$antihomomorphisms) with multiplicative derivations. More specifically, we study multiplicative (generalized)-derivations of prime rings that act as $n-$homomorphism or $n-$antiho\\momorphism.

\section{The Results}\label{Sec2}

We begin with the following observations in this subject, which we shall use frequently.


\begin{lemma}\label{lemma1}
Let $R$ be a prime ring and $I$ be a nonzero ideal of $R.$ Then for any $a,b\in R,$ $aIb=(0)$ implies $a=0$ or $b=0.$
\end{lemma}

\begin{lemma}\label{lemma2}
Let $R$ be a prime ring and $I$ be a nonzero ideal of $R.$ If for any fixed positive integer $n,$ $[x^{n},y^{n}]\in Z(R)$ for all $x,y\in I,$ then $R$ is commutative.
\begin{proof}
By hypothesis, we have $[[x^{n},y^{n}],r]=0$ for all $x,y\in I$ and $r\in R.$ It is well-known that $I$ and $R$ satisfy same polynomial identities. Thus, we have $[[x^{n},y^{n}],r]=0$ for all $x,y,r\in R.$ If possible suppose that $R$ is not commutative. By a famous result of Lanski \cite{Lanski1993}, $R\subseteq M_{n}(F),$ where $M_{n}(F)$ be a ring of $n\times n$ matrices, with $n\geq 2$ over a field $F.$ Moreover, $R$ and $M_{n}(F)$ satisfy the same polynomial identities. Choose $x=e_{11}, y=e_{12}+e_{22}$ and $r=e_{21},$ where $e_{ij}$ denotes matrix with $1$ at $ij-$entry and $0$ elsewhere. In this view, it follows that
\[
0=[[x^{n},y^{n}],r]=e_{11},
\]
a contradiction. Hence, $R$ is commutative.
\end{proof}
\end{lemma}

\begin{lemma}\label{lemma3}
Let $R$ be a ring and $\delta$ be a multiplicative derivation of $R.$ Then the
followings are true:
\begin{itemize}
\item[(i)] $\delta(0)=0.$
\item[(ii)] If $a\in Z(R),$ then $\delta(a)\in Z(R).$
\end{itemize}

\begin{proof}
$(i)$ $\delta(0)=\delta(0.0)=\delta(0).0+0.\delta(0)=0.$
$(ii)$ Let $a\in Z(R)$ and $\delta$ be a multiplicative derivation of $R$. Then for each $x\in R,$ we have
\[
\delta(ax)=\delta(a)x+a\delta(x)
\]
and
\[
\delta(ax)=\delta(xa)=\delta(x)a+x\delta(a).
\]
Together with above two equations, we get%
\[
\left[x,\delta(a)\right]=0~\text{for all }x\in R.
\]
Hence $\delta(a)\in R.$
\end{proof}
\end{lemma}


\begin{theorem}\label{theorem1}
  Let $R$ be a prime ring, $I$ a nonzero ideal of $R.$ Suppose that $F:R\to R$ is a multiplicative (generalized)-derivation associated with a multiplicative derivation $\delta$ of $R$ such that $F$ acts as $n-$homomorphism on $I.$ Then $\delta=0,$ and $F=0$ or there exists $\lambda\in C$ such that $F(x)=\lambda x$ for all $x\in R$ and $\lambda^{n-1}=1.$
\begin{proof}
By hypothesis, we have
\begin{equation}\label{Eq-1}
F(\prod_{i=1}^{n}x_{i})=\prod_{i=1}^{n}F(x_{i})
\end{equation}
for all $x_{i}\in I.$ On the other hand, we find
\begin{equation}\label{Eq-2}
F(\prod_{i=1}^{n}x_{i})=F(\prod_{i=1}^{n-1}x_{i})x_{n}+\prod_{i=1}^{n-1}x_{i}\delta(x_{n})
\end{equation}
for all $x_{i}\in I.$ Combining (\ref{Eq-1}) and (\ref{Eq-2}), we obtain
\begin{equation}\label{Eq-3}
  \prod_{i=1}^{n}F(x_{i})=F(\prod_{i=1}^{n-1}x_{i})x_{n}+\prod_{i=1}^{n-1}x_{i}\delta(x_{n})
\end{equation}
for all $x_{i}\in I.$ Replace $x_{n}$ by $x_{n}r$ in (\ref{Eq-3}), where $r\in R,$ we get
\[
  \prod_{i=1}^{n-1}F(x_{i})x_{n}\delta(r)=\prod_{i=1}^{n}x_{i}\delta(r).
\]
That is
\[
(\prod_{i=1}^{n-1}F(x_{i})-\prod_{i=1}^{n-1}x_{i})x_{n}\delta(r)=0.
\]
In view of Lemma \ref{lemma1}, we find that either $\prod_{i=1}^{n-1}F(x_{i})=\prod_{i=1}^{n-1}x_{i}$ or $\delta=0.$ Let us consider
\begin{equation}\label{Eq-4}
  \prod_{i=1}^{n-1}F(x_{i})=\prod_{i=1}^{n-1}x_{i}
\end{equation}
for all $x_{i}\in I.$ Replace $x_{n-1}$ by $x_{n-1}r$ in (\ref{Eq-4}), we find
\begin{equation}\label{Eq-5}
  \prod_{i=1}^{n-1}F(x_{i})r+\prod_{i=1}^{n-2}F(x_{i})x_{n-1}\delta(r)=\prod_{i=1}^{n-1}x_{i}r
\end{equation}
for all $x_{i}\in I$ and $r\in R.$ Right multiply (\ref{Eq-4}) by $r$ and subtract from (\ref{Eq-5}), we get
\[
\prod_{i=1}^{n-2}F(x_{i})x_{n-1}\delta(r)=0
\]
for all $x_{i}\in I$ and $r\in R.$ Again by invoking Lemma \ref{lemma1}, we find that either $\prod_{i=1}^{n-2}F(x_{i})=0$ or $\delta=0.$ But $\delta\neq 0,$ so we have $\prod_{i=1}^{n-2}F(x_{i})=0$ for all $x_{i}\in I.$ Substitute $x_{n-2}r$ in place of $x_{n-2}$ in above expression, where $r\in R,$ we find that $\prod_{i=1}^{n-3}F(x_{i})I\delta(r)=(0).$ By Lemma \ref{lemma1}, it follows that either $\prod_{i=1}^{n-3}F(x_{i})=0$ for all $x_{i}\in I$ or $\delta=0.$ But $\delta\neq 0,$ thus we have $\prod_{i=1}^{n-3}F(x_{i})=0$ for all $x_{i}\in I.$ Continuing in this way, we arrive at $F(x)=0$ for all $x\in I.$ Replace $x$ by $xr,$ where $r\in R,$ we get $x\delta(r)=0$ for all $x\in I$ and $r\in R.$ It implies that $\delta=0,$ which is a contradiction.
\par Let us now consider the latter case $\delta=0,$ we find that
\begin{equation}\label{Eq-6}
F(\prod_{i=1}^{n}x_{i})=F(x_{i})\prod_{i=2}^{n}x_{i}
\end{equation}
for all $x_{i}\in I.$ Combining (\ref{Eq-1}) and (\ref{Eq-6}), we obtain
\[
F(x_{1})(\prod_{i=2}^{n}F(x_{i})-\prod_{i=2}^{n}x_{i})=0
\]
for all $x_{i}\in I.$ Replace $x_{1}$ by $x_{1}r,$ where $r\in R,$ we may infer that
\[
F(x_{1})R(\prod_{i=2}^{n}F(x_{i})-\prod_{i=2}^{n}x_{i})=(0)
\]
 for all $x_{i}\in I.$ Since $R$ is prime, we find that either $F(x)=0$ for all $x\in I$ or $\prod_{i=2}^{n}F(x_{i})=\prod_{i=2}^{n}x_{i}$ for all $x_{i}\in I.$ It is straight forward to see that the former case implies $F=0.$ On the other side, we have
\begin{equation}\label{Eq-7}
  \prod_{i=2}^{n}F(x_{i})=\prod_{i=2}^{n}x_{i}
\end{equation}
for all $x_{i}\in I.$ Take $rx_{2}$ instead of $x_{2}$ in (\ref{Eq-7}), where $r\in R,$ we get
\begin{equation}\label{Eq-8}
  F(r)x_{2}\prod_{i=3}^{n}F(x_{i})=rx_{2}\prod_{i=3}^{n}x_{i}.
\end{equation}
Left multiply (\ref{Eq-7}) by $r$ and then subtract from (\ref{Eq-8}), we obtain
\[
(F(r)x_{2}-rF(x_{2}))\prod_{i=3}^{n}F(x_{i})=0
\]
for all $x_{i}\in I$ and $r\in R.$ Substitute $x_{2}s$ in place of $x_{2}$ in above equation, where $s\in R,$ we obtain
\[
(F(r)x_{2}-rF(x_{2}))R\prod_{i=3}^{n}F(x_{i})=(0)
\]
 for all $x_{i}\in I$ and $r\in R.$ It implies that either $F(r)x-rF(x)=0$ for all $x\in I$ and $r\in R$ or $\prod_{i=3}^{n}F(x_{i})=0$ for all $x_{i}\in I.$ One may observe that in both of these cases we get the situation $F(r)x-rF(x)=0$ for all $x\in I$ and $r\in R.$ Replace $x$ by $sx,$ we get $(F(r)s-rF(s))x=0$ for all $x\in I$ and $r,s\in R.$ By Lemma \ref{lemma1}, we get $F(r)s=rF(s)$ for all $r,s\in R.$ Replace $r$ by $rp,$ we get
 $F(r)p1_{R}(s)=1_{R}(r)pF(s)$ for all $r,s,p\in R,$ where $1_{R}$ is the identity mapping of $R.$ With the aid of a result of Bre$\check{s}$ar [\cite{Bresar1990}, Lemma], it follows that there exists some $\lambda\in C$ such that $F=\lambda 1_{R}$ and hence $F(x)=\lambda x$ for all $x\in R.$ In view of our hypothesis, we have $\lambda\prod_{i=1}^{n}x_{i}=\prod_{i=1}^{n}\lambda x_{i}.$ It forces that $\lambda^{n-1}=1.$ It completes the proof.
\end{proof}
\end{theorem}

\begin{corollary}\label{corollary1}\textsc{[\cite{Gusic2005}, Theorem 1($a$)]}
Let $R$ be an associative prime ring, $I$ a nonzero ideal of $R.$ Suppose that $F:R\to R$ is a multiplicative (generalized)-derivation associated with a multiplicative derivation $\delta$ of $R$ such that $F$ acts a homomorphism on $I.$ Then $\delta=0,$ and $F=0$ or $F(x)=x$ for all $x\in R.$
\end{corollary}

In spirit of a result of Gusi\'{c} [\cite{Gusic2005}, Theorem 1($b$)], it is natural to investigate multiplicative (generalized)-derivations that act as $n-$antihomomorphisms. However, we could not get this result in its complete form, but we obtain the following:

\begin{theorem}\label{theorem2}
  Let $R$ be a prime ring, $I$ a nonzero ideal of $R.$ Suppose that $F:R\to R$ is a multiplicative (generalized)-derivation associated with a multiplicative derivation $\delta$ of $R$ such that $F$ acts as $n-$antihomomorphism on $I.$ If $F=\delta,$ then $\delta(x)^{n-1}\in Z(R)$ for all $x\in I.$ Moreover, if $\delta$ is additive, then either $\delta=0$ or $R$ is commutative or $R$ is an order in a $4-$dimensional simple algebra.
\begin{proof}
By hypothesis, we have
\begin{equation}\label{Eq-9}
  F(\prod_{i=1}^{n})=F(x_{n})F(x_{n-1})\cdots F(x_{2})F(x_{1})
\end{equation}
for all $x_{i}\in I.$ On the other hand, we may infer that
\begin{equation}\label{Eq-10}
  F(\prod_{i=1}^{n})=F(x_{1})\prod_{i=2}^{n}x_{i}+\sum_{i=2}^{n}(\prod_{j=1}^{i-1}x_{j}\delta(x_{i})\prod_{k=i+1}^{n}x_{k})
\end{equation}
for all $x_{i}\in I.$ Combining (\ref{Eq-9}) and (\ref{Eq-10}), we find that
\begin{equation}\label{Eq-11}
  F(x_{n})\cdots F(x_{1})=F(x_{1})\prod_{i=2}^{n}x_{i}+\sum_{i=2}^{n}(\prod_{j=1}^{i-1}x_{j}\delta(x_{i})\prod_{k=i+1}^{n}x_{k})
\end{equation}
for all $x_{i}\in I.$ Replace $x_{1}$ by $x_{1}x_{n}$ in (\ref{Eq-11}), we obtain
\begin{equation}\label{Eq-12}
\begin{split}
  F(x_{n})\cdots F(x_{2})F(x_{1})x_{n}+F(x_{n})\cdots F(x_{2})x_{1}\delta(x_{n})
  =F(x_{1})x_{n}\prod_{i=2}^{n}x_{i}
  \\+x_{1}\delta(x_{n})\prod_{i=2}^{n}x_{i}
  +x_{1}x_{n}\delta(x_{2})\prod_{i=3}^{n}x_{i}
  +x_{1}x_{n}\sum_{i=3}^{n}(\prod_{j=2}^{i-1}
  x_{j}\delta(x_{i})\prod_{k=i+1}^{n}x_{k})
\end{split}
\end{equation}
for all $x_{i}\in I.$ Using (\ref{Eq-9}) in (\ref{Eq-12}), we get
\[
\begin{split}
  F(\prod_{i=1}^{n}x_{i})x_{n}+F(x_{n})\cdots F(x_{2})x_{1}\delta(x_{n})=F(x_{1})x_{n}\prod_{i=2}^{n}x_{i}+x_{1}
  \delta(x_{n})\prod_{i=2}^{n}x_{i}
  \\+x_{1}x_{n}\delta(x_{2})\prod_{i=3}^{n}x_{i}+x_{1}x_{n}\sum_{i=3}^{n}
  (\prod_{j=2}^{i-1}x_{j}\delta(x_{i})\prod_{k=i+1}^{n}x_{k})
\end{split}
\]
for all $x_{i}\in I.$ It implies that
\[
\begin{split}
(F(x_{1})\prod_{i=2}^{n}x_{i}+\sum_{i=2}^{n}(\prod_{j=1}^{i-1}x_{j}\delta(x_{i})\prod_{k=i+1}^{n}x_{k}))x_{n}+F(x_{n})\cdots F(x_{2})x_{1}\delta(x_{n})
\\=F(x_{1})x_{n}\prod_{i=2}^{n}x_{i}+x_{1}\delta(x_{n})\prod_{i=2}^{n}x_{i}+x_{1}x_{n}\delta(x_{2})\prod_{i=3}^{n}x_{i}+x_{1}x_{n}
\\ \sum_{i=3}^{n}(\prod_{j=2}^{i-1}x_{j}\delta(x_{i})\prod_{k=i+1}^{n}x_{k})
\end{split}
\]
for all $x_{i}\in I.$ In particular, for $x_{1}=x$ and $x_{2}=x_{3}=\cdots=x_{n}=y,$ we find
\[
\begin{split}
F(x)y^{n}+x(\sum_{i=0}^{n-2}y^{i}\delta(y)y^{n-1-i})+F(y)^{n-1}x\delta(y)=F(x)y^{n}+x\delta(y)y^{n-1}
\\+xy\delta(y)y^{n-2}+x(\sum_{i=2}^{n-1}y^{i}\delta(y)y^{n-1-i})
\end{split}
\]
for all $x,y\in I.$ It yields that
\begin{equation}\label{Eq-13}
F(y)^{n-1}x\delta(y)=xy^{n-1}\delta(y)
\end{equation}
for all $x,y\in I.$ Replace $x$ by $rx,$ where $r\in R$ in (\ref{Eq-13}), we get
\begin{equation}\label{Eq-14}
F(y)^{n-1}rx\delta(y)=rxy^{n-1}\delta(y).
\end{equation}
Left multiply (\ref{Eq-13}) by $r$ and combine with (\ref{Eq-14}), we obtain $[F(y)^{n-1},r]x\delta(y)=0$ for all $x,y\in I$ and $r\in R.$
\par In particular, we take $F=\delta.$ Thus we have $[\delta(y)^{n-1},r]x\delta(y)=0$ for all $x,y\in I$ and $r\in R.$ Since $R$ is a prime ring,
it follows that for each $y\in I,$ either $[\delta(y)^{n-1},r]=0$ for all $r\in R$ or $\delta(y)=0.$ In each case we have $[\delta(y)^{n-1},r]=0$
for all $y\in I$ and $r\in R,$ i.e., $\delta(y)^{n-1}\in Z(R)$ for all $y\in I.$ In case $\delta$ is additive, we are done by
[\cite{Chang2009}, Theorem B].
\end{proof}
\end{theorem}

\begin{corollary}\label{corollary2}\textsc{[\cite{Gusic2005}, Theorem 1($b$)]}
Let $R$ be an associative prime ring, $I$ a nonzero ideal of $R.$ Suppose that $F:R\to R$ is a multiplicative (generalized)-derivation associated with a multiplicative derivation $\delta$ of $R$ such that $F$ acts a homomorphism on $I.$ Then $\delta=0,$ and $F=0$ or $F(x)=x$ for all $x\in R.$
\begin{proof}
  For $n=2,$ in view of equation (\ref{Eq-13}) and (\ref{Eq-14}), we have $[F(y),t]x\delta(y)=0$ for all $x,y,t\in I.$ This same expression appeared in the beginning of the proof of Theorem 1(b) in \cite{Gusic2005}, hence the conclusion follows in the same way.
\end{proof}
\end{corollary}

\begin{definition}
Let $F:R\to R$ be a function. Then $F$ is called right multiplicative (generalized)-derivation of $R$ if it satisfies
\[
F(xy)=F(x)y+x\delta(y)
\]
for all $x,y\in R$ and $\delta$ is any mapping of $R.$ And $F$ is called left multiplicative (generalized)-derivation of $R$ if it satisfies
\[
F(xy)=\delta(x)y+xF(y)
\]
for all $x,y\in R$ and $\delta$ is any mapping of $R.$ Then it is not difficult to see that the associated mapping $\delta$ of right and left multiplicative (generalized)-derivation $F$ is a multiplicative derivation. Now, $F$ is said to be two-sided multiplicative (generalized)-derivation of $R$ if it satisfies
\begin{eqnarray*}
  F(xy) &=& F(x)y+x\delta(y) \\
   &=& \delta(x)y+xF(y)
\end{eqnarray*}
for all $x,y\in R,$ where $\delta$ is a multiplicative derivation of $R.$
\end{definition}

\begin{theorem}\label{theorem3}
  Let $R$ be a prime ring, $I$ a nonzero ideal of $R.$ Suppose that $F:R\to R$ is a two-sided multiplicative (generalized)-derivation associated with a multiplicative derivation $\delta$ of $R$ such that $F$ acts as $n-$antihomomorphism on $I.$ Then $\delta=0,$ and $F=0$ or there exists $\lambda\in C$ such that $F(x)=\lambda x$ for all $x\in R$ and $\lambda^{n-1}=1$ (in this case $R$ should be commutative).
\begin{proof}
From equation (\ref{Eq-13}), we have $F(y)^{n-1}x\delta(y)=xy^{n-1}\delta(y)$ for all $x,y\in I.$ Take $F(z)x$ in place of $x$ in this equation, we get
\begin{eqnarray*}
  F(y)^{n-1}F(z)x\delta(y) &=& F(z)xy^{n-1}\delta(y)\\
  F(zy^{n-1})x\delta(y) &=& F(z)xy^{n-1}\delta(y) \\
  F(z)y^{n-1}x\delta(y)+z\delta(y^{n-1})x\delta(y) &=& F(z)xy^{n-1}\delta(y)
\end{eqnarray*}
for all $x,y,z\in I.$ It implies that
\begin{equation}\label{Eq-15}
  F(z)[y^{n-1},x]\delta(y)+z\delta(y^{n-1})x\delta(y)=0
\end{equation}
for all $x,y,z\in I.$ Replace $z$ by $rz$ in (\ref{Eq-15}), where $r\in R,$ we get
\[
\delta(r)z[y^{n-1},x]\delta(y)+rF(z)[y^{n-1},x]\delta(y)+rz\delta(y^{n-1})x\delta(y)=0.
\]
Using (\ref{Eq-15}), we find $\delta(r)z[y^{n-1},x]\delta(y)=0$ for all $x,y,z\in I$ and $r\in R.$ In view of Lemma \ref{lemma1}, it implies that either $\delta=0$ or $[y^{n-1},x]\delta(y)=0$ for all $x,y\in I.$ Assume that $[y^{n-1},x]\delta(y)=0$ for all $x,y\in I.$ It implies that for each $y\in I,$ either $y^{n-1}\in Z(R)$ or $\delta(y)=0.$ Together these both cases (using Lemma \ref{lemma3}) imply that $\delta(y^{n-1})\in Z(R)$ for all $y\in I.$
\par We now consider
\begin{eqnarray*}
  F(xy^{n-1}) &=& F(x)y^{n-1}+x\delta(y^{n-1}) \\
  F(xy^{n-1}) &=& F(y)^{n-1}F(x)
\end{eqnarray*}
for all $x,y\in I.$ Thus we have
\begin{eqnarray}\label{Eq-16}
  F(y)^{n-1}F(x) &=& F(x)y^{n-1}+x\delta(y^{n-1})\notag \\
   &=& F(x)y^{n-1}+\delta(y^{n-1})x.
\end{eqnarray}
Take $xz$ in place of $x$ in (\ref{Eq-16}), we find
\begin{equation}\label{Eq-17}
  F(y)^{n-1}F(x)z+F(y)^{n-1}x\delta(z)=F(x)zy^{n-1}+x\delta(z)y^{n-1}+\delta(y^{n-1})xz
\end{equation}
for all $x,y,z\in I.$ Using (\ref{Eq-16}), it implies that
\begin{equation}\label{Eq-18}
F(y)^{n-1}x\delta(z)=F(x)[z,y^{n-1}]+x\delta(z)y^{n-1}
\end{equation}
for all $x,y,z\in I.$ Replace $x$ by $rx$ in (\ref{Eq-18}), where $r\in R,$ we get
\[
F(y)^{n-1}rx\delta(z)=rF(x)[z,y^{n-1}]+\delta(r)x[z,y^{n-1}]+rx\delta(z)y^{n-1}.
\]
Using (\ref{Eq-18}), we have
\begin{equation}\label{Eq-19}
[F(y)^{n-1},r]x\delta(z)=\delta(r)x[z,y^{n-1}]
\end{equation}
for all $x,y,z\in I$ and $r\in R.$ Replace $z$ by $zw^{n-1}$ in (\ref{Eq-19}), we get
\[
\begin{split}
[F(y)^{n-1},r]x\delta(z)w^{n-1}+[F(y)^{n-1},r]xz\delta(w^{n-1})=\delta(r)x[z,y^{n-1}]w^{n-1}
\\+\delta(r)xz[y^{n-1},w^{n-1}]
\end{split}
\]
for all $x,y,z,w\in I$ and $r\in R.$ Equation (\ref{Eq-19}) reduces it to
\begin{equation}\label{Eq-20}
  \delta(w^{n-1})[F(y)^{n-1},r]xz=\delta(r)xz[y^{n-1},w^{n-1}]
\end{equation}
for all $x,y,z,w\in I$ and $r\in R.$
Take $zs$ in place of $z$ in (\ref{Eq-20}), where $s\in R,$ we find
\[
  \delta(w^{n-1})[F(y)^{n-1},r]xzs=\delta(r)xzs[y^{n-1},w^{n-1}]
\]
for all $x,y,z,w\in I$ and $r,s\in R.$ Using (\ref{Eq-20}) in the above expression, we obtain $\delta(r)xz[[w^{n-1},\\y^{n-1}],s]=0$ for all $x,y,z,w\in I$ and $r,s\in R.$ It forces that either $\delta=0$ or $[w^{n-1},y^{n-1}]\in Z(R)$ for all $y,w\in I.$ But $\delta\neq 0,$ thus we have $[w^{n-1},y^{n-1}]\in Z(R)$ for all $y,w\in I.$ In view of Lemma \ref{lemma2}, $R$ is commutative. Therefore, $F$ is just $n-$homomorphism of $R$ and hence by Theorem \ref{theorem1}, we get $\delta=0,$ a contradiction.
\par On the other hand, we assume that $\delta=0.$ Relation (\ref{Eq-10}) implies that
\[
F(x_{1}x_{2}\cdots x_{n})=F(x_{1})x_{2}\cdots x_{n}
\]
for all $x_{i}\in I.$ Using this relation, we obtain
\begin{eqnarray*}
  F(x_{1})x_{2}x_{3}\cdots x_{n-1}x_{n}x_{n+1} &=& F(x_{1}x_{2}\cdots x_{n-1}x_{n})x_{n+1} \\
   &=& F(x_{n})F({x_{n-1}})\cdots F(x_{2})F(x_{1})x_{n+1} \\
   &=& F(x_{n})F({x_{n-1}})\cdots F(x_{2})F(x_{1}x_{n+1}) \\
   &=& F(x_{1}x_{n+1}x_{2}\cdots x_{n}) \\
   &=& F(x_{1})x_{n+1}x_{2}\cdots x_{n}
\end{eqnarray*}
for all $x_{i}\in I.$ It gives
\[
F(x_{1})[x_{2}\cdots x_{n},x_{n+1}]=0
\]
for all $x_{i}\in I.$ Thus we have either $F(x)=0$ for all $x\in I$ or $[x_{2}\cdots x_{n},x_{n+1}]=0$ for all $x_{i}\in I.$ The first case implies $F=0.$ In the latter case we find that $R$ is commutative and hence $F$ acts as $n-$homomorphism on $I.$ We are done by Theorem \ref{theorem1}.

\end{proof}
\end{theorem}

\section*{Acknowledgement}
I would like to thank Prof. Ne\c{s}et AYDIN for reading the earlier draft of the manuscript and suggesting Lemma \ref{lemma3}. 

\end{document}